\documentclass[10pt]{article}
\usepackage{amsmath,amssymb,amsfonts,amsthm,mathtools,epsfig}
\usepackage{tikz}
\usepackage{caption,subcaption,verbatim}
\usepackage{setspace}
\usepackage{xurl}
\usepackage{hyperref}
\hypersetup{hidelinks}
\usepackage[normalem]{ulem}

\usepackage{cancel}

\usepackage[bottom=1.5in]{geometry}

\newtheorem{thm}{Theorem}
\newtheorem{theorem}[thm]{Theorem}
\newtheorem{lemma}[thm]{Lemma}

\newtheorem{question}[thm]{Question}
\newtheorem{obs}[thm]{Observation}
\theoremstyle{definition}

\title{\textbf{An infinite family of bipartite hypotraceable graphs}}

\author{
Jorik Jooken
\thanks{Mathematical Institute, Leiden University, 2333 CA Leiden, 
the Netherlands.\\ E-mail: 
\protect\href{mailto:j.jooken@math.leidenuniv.nl}
{\protect\mbox{\protect\nolinkurl{j.jooken@math.leidenuniv.nl}}}}
\thanks{Department of Computer Science, KU Leuven Kulak, 8500 Kortrijk, 
Belgium.\\ E-mail: 
\protect\href{mailto:jorik.jooken@kuleuven.be}
{\protect\mbox{\protect\nolinkurl{jorik.jooken@kuleuven.be}}}}
}

\date{}

\begin{document}

	\maketitle

	\begin{abstract}
	A graph is hypotraceable if it has no Hamiltonian path but every vertex-deleted subgraph does. In 1978, Gr\"otschel asked whether a bipartite hypotraceable graph exists. We resolve this long-standing question in a strong sense by proving that such graphs exist for every sufficiently large even order. Moreover, we show that bipartite hypotraceable graphs can have $\left(1/4+o(1)\right)n^2$ edges, asymptotically attaining the maximum possible number of edges of a bipartite graph, and that there are at least $2^{(1/4+o(1))n^2}$ nonisomorphic examples of order $n$, as $n$ tends to infinity through even integers. Finally, we prove that every bipartite graph occurs as an induced subgraph of a bipartite hypotraceable graph.
			\vskip 3mm

			\end{abstract}

	\section{Introduction}\label{sec:intro}

Hypohamiltonian and hypotraceable graphs have been studied for more than half a century~\cite{Sousselier63,Thomassen74}. A graph is \textit{Hamiltonian} if it contains a spanning cycle, referred to as a \textit{Hamiltonian cycle}. Similarly,
a graph is \textit{traceable} if it contains a spanning path, referred to as a \textit{Hamiltonian path}. A graph $G$ is called \textit{hypohamiltonian} if $G$ is not Hamiltonian, but $G-v$ is Hamiltonian for every vertex $v \in V(G)$. In the current paper, we shall mostly concern ourselves with the related notion of \textit{hypotraceable graphs}: $G$ is hypotraceable if it is not traceable, but $G-v$ is traceable for every vertex $v \in V(G)$.

Initially, the existence of hypotraceable graphs was even conjectured to be impossible, until Horton found an example on $40$ vertices in 1973~\cite{Wiener18}. Thomassen~\cite{Thomassen74} subsequently proposed a construction to obtain hypotraceable graphs starting from hypohamiltonian graphs, thereby obtaining a smaller example on $34$ vertices. He later constructed infinite families of planar and planar cubic examples~\cite{Thomassen76,Thomassen81}. Such constructions were further developed in subsequent literature, including the planar examples of Jooyandeh, McKay, \"Osterg\aa rd, Pettersson and Zamfirescu~\cite{JMOPZ17} and Wiener's extension using almost hypohamiltonian graphs~\cite{Wiener18}. These graphs also play a pivotal role in combinatorial optimisation. More precisely, Gr\"otschel~\cite{Groetschel80} showed that certain classes of hypohamiltonian and hypotraceable graphs lead to facet-defining inequalities for the monotone symmetric travelling salesman polytope.

Despite significant attention in the literature on constructing such graphs, the bipartite case has resisted a basic existence question. It is easy to see that a bipartite hypohamiltonian graph cannot exist, because a Hamiltonian cycle in a bipartite graph must use equally many vertices from its two partite sets. For Hamiltonian paths on the other hand, the sizes may differ by one. This leaves the possibility of bipartite hypotraceable graphs open, although their two partite sets must have equal size and their order must therefore be even. In 1978, Gr\"otschel famously asked whether such graphs exist\footnote{To avoid the trivial solution of the $2$-vertex graph consisting of two isolated vertices, we will require $|V(G)| \geq 3$.}.
\begin{question}[Gr\"otschel, 1978]
Does there exist a bipartite hypotraceable graph?
\end{question}
The question appears in Bondy's chapter in the \emph{Handbook of Combinatorics}~\cite[Problem~4.56]{Bondy95}. It also appears in Bondy and Murty's book \emph{Graph Theory}~\cite[Appendix~A, Problem~92]{BondyMurty08}, where it is actually formulated as the conjecture that no such graph exists. This question has continued to motivate work on related spanning structures: Goedgebeur, Renders and Zamfirescu~\cite{GRZ24} recall this question in their work on bipartite $K_2$-hypohamiltonian graphs. Goedgebeur, Renders, Wiener and Zamfirescu~\cite{GRWZ26} relaxed it in their work on minimum leaf spanning trees and fault-tolerant networks, while Wiener, Yokota and Zamfirescu~\cite[Problem~4]{WYZ21} proposed another relaxation of Gr\"otschel's question in terms of spanning spiders.

In the current paper, we answer Gr\"otschel's question affirmatively and thereby disprove the nonexistence conjecture from Bondy and Murty's book. More precisely, we prove that a bipartite hypotraceable graph exists for every even order $n \geq 200$. Our examples cover a wide range of sizes, from graphs having $3n/2+10$ edges to graphs having $\left(1/4+o(1)\right)n^2$ edges, asymptotically attaining the maximum possible number of edges of a bipartite graph. We also prove that there are at least $2^{(1/4+o(1))n^2}$ nonisomorphic bipartite hypotraceable graphs of order $n$, as $n$ tends to infinity through even integers, and that every bipartite graph occurs as an induced subgraph of some bipartite hypotraceable graph.

    \section{Results}\label{sec:results}

\subsection{An initial building block}\label{subsec:D}

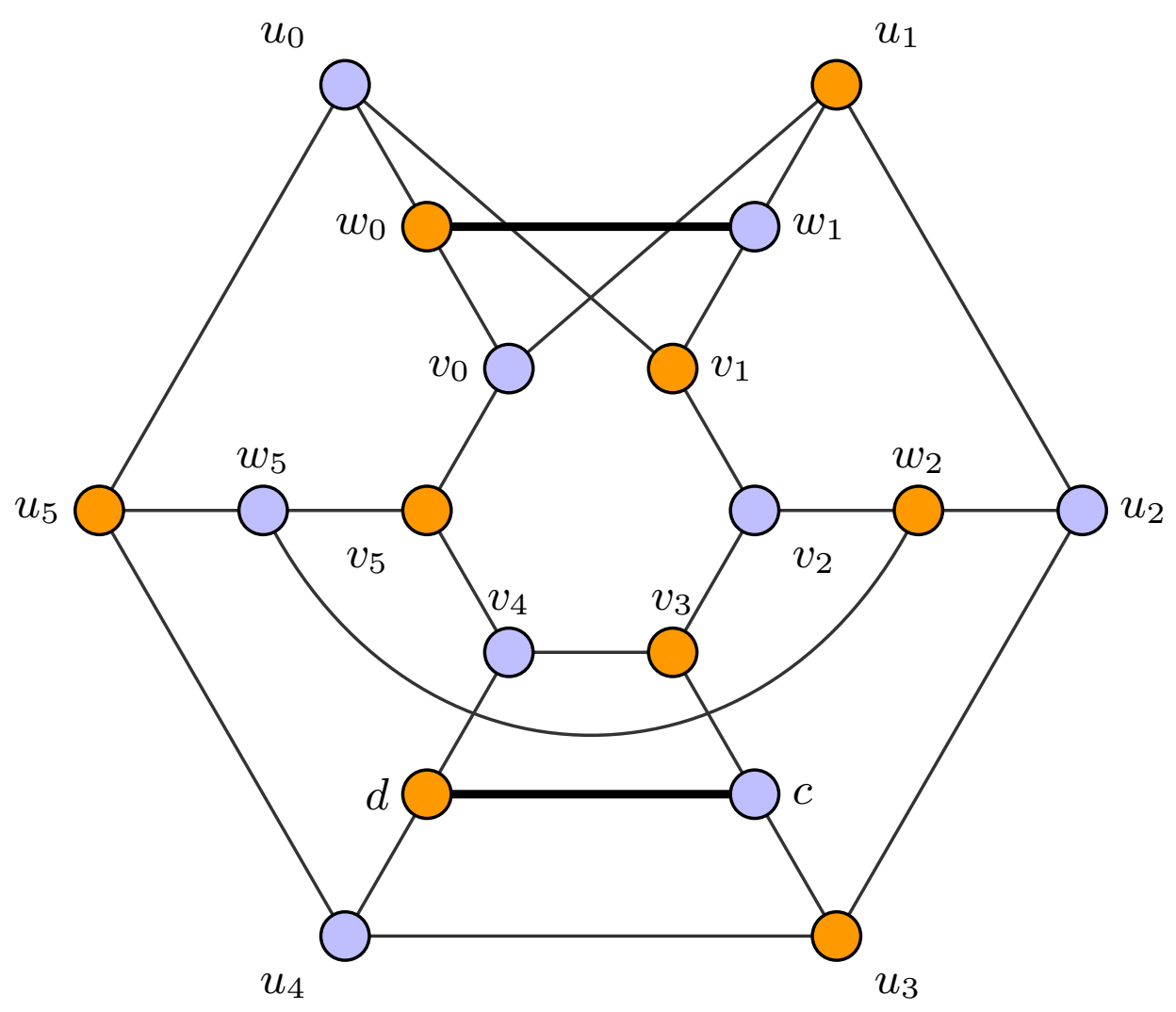
\begin{figure}[htbp]
\centering
\begin{tikzpicture}[
  x=1.15cm,y=1.15cm,
  edge/.style={draw=black!80,line width=0.65pt},
  special/.style={draw=black,line width=1.7pt},
  vertex/.style={circle,draw=black,line width=0.65pt,inner sep=0pt,minimum size=3.4mm},
  first/.style={vertex,fill=blue!25},
  second/.style={vertex,fill=orange!80!yellow}
]
\coordinate (u0) at (-1.5,2.598076);
\coordinate (u1) at (1.5,2.598076);
\coordinate (u2) at (3,0);
\coordinate (u3) at (1.5,-2.598076);
\coordinate (u4) at (-1.5,-2.598076);
\coordinate (u5) at (-3,0);
\coordinate (w0) at (-1,1.732051);
\coordinate (w1) at (1,1.732051);
\coordinate (w2) at (2,0);
\coordinate (c) at (1,-1.732051);
\coordinate (d) at (-1,-1.732051);
\coordinate (w5) at (-2,0);
\coordinate (v0) at (-0.5,0.866025);
\coordinate (v1) at (0.5,0.866025);
\coordinate (v2) at (1,0);
\coordinate (v3) at (0.5,-0.866025);
\coordinate (v4) at (-0.5,-0.866025);
\coordinate (v5) at (-1,0);
\foreach \u/\v in {w1/v1,w1/u1,v2/v1,v2/v3,v2/w2,u0/v1,u0/u5,u0/w0,u2/u1,u2/u3,u2/w2,v0/u1,v0/v5,v0/w0,c/v3,c/u3,v4/v3,v4/d,v4/v5,u4/u3,u4/d,u4/u5,w5/v5,w5/u5}
  \draw[edge] (\u) -- (\v);
\draw[edge] (w5) .. controls (-1.1,-1.83) and (1.1,-1.83) .. (w2);
\draw[special] (w0) -- (w1) (d) -- (c);
\node[first] at (u0) {};
\node[second] at (u1) {};
\node[first] at (u2) {};
\node[second] at (u3) {};
\node[first] at (u4) {};
\node[second] at (u5) {};
\node[second] at (w0) {};
\node[first] at (w1) {};
\node[second] at (w2) {};
\node[first] at (c) {};
\node[second] at (d) {};
\node[first] at (w5) {};
\node[first] at (v0) {};
\node[second] at (v1) {};
\node[first] at (v2) {};
\node[second] at (v3) {};
\node[first] at (v4) {};
\node[second] at (v5) {};
\node[above left=4pt,font=\small] at (u0) {$u_0$};
\node[above right=4pt,font=\small] at (u1) {$u_1$};
\node[right=4pt,font=\small] at (u2) {$u_2$};
\node[below right=4pt,font=\small] at (u3) {$u_3$};
\node[below left=4pt,font=\small] at (u4) {$u_4$};
\node[left=4pt,font=\small] at (u5) {$u_5$};
\node[left=4pt,font=\small] at (w0) {$w_0$};
\node[right=4pt,font=\small] at (w1) {$w_1$};
\node[above=4pt,font=\small] at (w2) {$w_2$};
\node[right=4pt,font=\small] at (c) {$c$};
\node[left=4pt,font=\small] at (d) {$d$};
\node[above=4pt,font=\small] at (w5) {$w_5$};
\node[left=4pt,font=\small] at (v0) {$v_0$};
\node[right=4pt,font=\small] at (v1) {$v_1$};
\node[below right=4pt,font=\small] at (v2) {$v_2$};
\node[above=4pt,font=\small] at (v3) {$v_3$};
\node[above=4pt,font=\small] at (v4) {$v_4$};
\node[below left=4pt,font=\small] at (v5) {$v_5$};
\end{tikzpicture}
\caption{The $18$-vertex cubic bipartite graph $B$. No Hamiltonian cycle contains simultaneously the edges $w_0w_1$ and $cd$ (shown in bold).}
\label{fig:B}
\end{figure}

Let $B$ be the cubic bipartite $18$-vertex graph shown in Figure~\ref{fig:B}. This graph appears in the paper by Ellingham and Horton~\cite{EH83}, who attribute its use to an earlier construction by Ellingham~\cite{Ellingham82}. This graph has the following useful property.

\begin{obs}[Ellingham and Horton~{\cite{EH83}}]\label{obs:B}
No Hamiltonian cycle of $B$ contains both $w_0w_1$ and $cd$.
\end{obs}

Let $D$ be the graph that is obtained from $B$ by subdividing the edge $w_0w_1$ two times, thereby replacing it with the path $(w_0,a,b,w_1)$. The vertices $a, b, c$ and $d$ from the graph $D$ will play a special role and we will refer to these vertices as \textit{terminals}. The graph $D$ is a bipartite graph with partite sets $\mathcal A=\{a,c,u_0,u_2,u_4,v_0,v_2,v_4,w_1,w_5\}$ and $\mathcal B=\{b,d,u_1,u_3,u_5,v_1,v_3,v_5,w_0,w_2\}$. We now record several structural properties of $D$.

\begin{lemma}\label{lem:D}
The following statements hold for the graph $D$:
\begin{enumerate}
\item[(i)] $D$ contains a Hamiltonian $ab$-path, a Hamiltonian $ad$-path, a Hamiltonian $cb$-path, but no Hamiltonian $cd$-path.
\item[(ii)] For each two vertex-disjoint paths that together contain all vertices in $D$ and have $a,b,c$ and $d$ as their four endpoints, one must be an $ab$-path and the other one must be a $cd$-path.
\item[(iii)] For every $z\in\mathcal A$, the graph $D-z$ has a Hamiltonian $bd$-path, and for every $z\in\mathcal B$, the graph $D-z$ has a Hamiltonian $ac$-path.
\end{enumerate}
\end{lemma}

\begin{proof}
We first prove (i). The following three paths are respectively Hamiltonian $ab$-paths, Hamiltonian $ad$-paths and Hamiltonian $cb$-paths: 
\[
\begin{aligned}
&(a,w_0,u_0,v_1,v_2,v_3,c,u_3,u_2,w_2,w_5,u_5,u_4,d,v_4,v_5,v_0,u_1,w_1,b),\\
&(a,b,w_1,v_1,v_2,v_3,c,u_3,u_4,u_5,u_0,w_0,v_0,u_1,u_2,w_2,w_5,v_5,v_4,d),\\
&(c,v_3,v_2,v_1,w_1,u_1,v_0,v_5,v_4,d,u_4,u_3,u_2,w_2,w_5,u_5,u_0,w_0,a,b).
\end{aligned}
\]
Suppose for the sake of obtaining a contradiction that $Q$ is a Hamiltonian $cd$-path in $D$. Since $a$ and $b$ have degree two, this Hamiltonian path must contain the subpath $(w_0,a,b,w_1)$. But then one obtains a Hamiltonian cycle of $B$ by adding the edge $cd$ and suppressing $a$ and $b$, which leads to a contradiction by Observation~\ref{obs:B} since this cycle contains both $w_0w_1$ and $cd$.

Next, we prove (ii). If an $ac$-path and a $bd$-path are vertex-disjoint and together contain all vertices of $D$, then one obtains a Hamiltonian cycle of $B$ by adding the edges $ab$ and $cd$ and suppressing $a$ and $b$. We obtain the same conclusion if there is such an $ad$-path and a $bc$-path. In both cases, we again obtain a contradiction by Observation~\ref{obs:B}.

Finally, we prove (iii). By symmetry of $D$, it suffices to prove the first statement. The required Hamiltonian $bd$-paths are shown in Table~\ref{tab:deletedD}.

\begin{center}
\begin{minipage}{\linewidth}
\centering
\small
\renewcommand{\arraystretch}{1.25}
\begin{tabular}{@{}cl@{}}
\hline
$z$ & Hamiltonian $bd$-path in $D-z$\\
\hline
$a$ & $(b,w_1,v_1,v_2,v_3,c,u_3,u_4,u_5,u_0,w_0,v_0,u_1,u_2,w_2,w_5,v_5,v_4,d)$\\
$c$ & $(b,a,w_0,u_0,v_1,w_1,u_1,v_0,v_5,w_5,u_5,u_4,u_3,u_2,w_2,v_2,v_3,v_4,d)$\\
$u_0$ & $(b,a,w_0,v_0,u_1,w_1,v_1,v_2,w_2,u_2,u_3,c,v_3,v_4,v_5,w_5,u_5,u_4,d)$\\
$u_2$ & $(b,a,w_0,u_0,v_1,w_1,u_1,v_0,v_5,v_4,v_3,v_2,w_2,w_5,u_5,u_4,u_3,c,d)$\\
$u_4$ & $(b,a,w_0,u_0,u_5,w_5,w_2,u_2,u_3,c,v_3,v_2,v_1,w_1,u_1,v_0,v_5,v_4,d)$\\
$v_0$ & $(b,a,w_0,u_0,v_1,w_1,u_1,u_2,w_2,v_2,v_3,c,u_3,u_4,u_5,w_5,v_5,v_4,d)$\\
$v_2$ & $(b,a,w_0,u_0,v_1,w_1,u_1,v_0,v_5,v_4,v_3,c,u_3,u_2,w_2,w_5,u_5,u_4,d)$\\
$v_4$ & $(b,a,w_0,u_0,v_1,w_1,u_1,v_0,v_5,w_5,u_5,u_4,u_3,u_2,w_2,v_2,v_3,c,d)$\\
$w_1$ & $(b,a,w_0,u_0,v_1,v_2,v_3,c,u_3,u_4,u_5,w_5,w_2,u_2,u_1,v_0,v_5,v_4,d)$\\
$w_5$ & $(b,a,w_0,u_0,u_5,u_4,u_3,u_2,w_2,v_2,v_1,w_1,u_1,v_0,v_5,v_4,v_3,c,d)$\\
\hline
\end{tabular}
\captionof{table}{A Hamiltonian $bd$-path in the graph $D-z$ for every vertex $z \in \mathcal A$.}
\label{tab:deletedD}
\end{minipage}
\end{center}

\end{proof}

\subsection{Combining building blocks}\label{subsec:X}

Next, we will combine two copies of $D$ and a bipartite graph $S$ into a new graph having four terminals. This graph $S$ offers us considerable freedom in our construction, from which we will later derive several interesting consequences. 

Let $D_L$ and $D_R$ be two copies of $D$. For every $i\in\{L,R\}$, we denote the terminals of $D_i$ by $a_i,b_i,c_i,d_i$ and its partite sets by $\mathcal A_i$ and $\mathcal B_i$. Let $S$ be a bipartite graph with partite sets $U$ and $W$ such that $|U|=|W|$. We allow $S$ to have no vertices, in which case $U=W=\varnothing$. If $S$ is nonempty, choose vertices $x_L,x_R\in U$ and $y_L,y_R\in W$ satisfying the following conditions:
\begin{enumerate}
\renewcommand{\theenumi}{\alph{enumi}}
\renewcommand{\labelenumi}{(\theenumi)}
\item\label{cond:S-paths} $S$ has a Hamiltonian $x_Ly_R$-path and a Hamiltonian $y_Lx_R$-path.
\item\label{cond:S-deletion} For every $z\in U$, the graph $S-z$ has a Hamiltonian $y_Ly_R$-path and for every $z\in W$, the graph $S-z$ has a Hamiltonian $x_Lx_R$-path.
\end{enumerate}

When $S$ is nonempty, we require the four chosen vertices $x_L,x_R,y_L$ and $y_R$ to be distinct, with one exception: $S$ is also allowed to be isomorphic to $K_2$, in which case $x_L=x_R$ and $y_L=y_R$ and the paths in~(\ref{cond:S-deletion}) each consist of a single vertex.

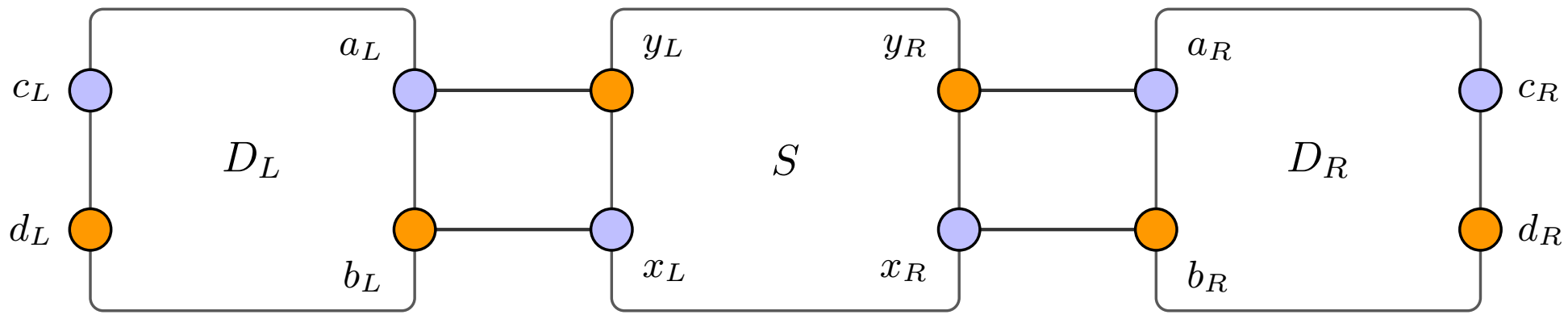
\begin{figure}[htbp]
\centering
\begin{tikzpicture}[
 x=0.95cm,y=0.95cm,
 block/.style={draw=black!65,line width=0.65pt,rounded corners=3pt},
 edge/.style={draw=black!80,line width=0.8pt},
 vertex/.style={circle,draw=black,line width=0.65pt,inner sep=0pt,minimum size=3.4mm},
 first/.style={vertex,fill=blue!25},
 second/.style={vertex,fill=orange!80!yellow},
 every label/.style={font=\small}
]
\draw[block] (0,0) rectangle (2.8,2.6);
\draw[block] (4.5,0) rectangle (7.5,2.6);
\draw[block] (9.2,0) rectangle (12,2.6);
\node at (1.4,1.3) {$D_L$};
\node at (6,1.3) {$S$};
\node at (10.6,1.3) {$D_R$};
\draw[edge] (2.8,1.9)--(4.5,1.9) (2.8,0.7)--(4.5,0.7);
\draw[edge] (7.5,1.9)--(9.2,1.9) (7.5,0.7)--(9.2,0.7);
\node[first,label=left:$c_L$] at (0,1.9) {};
\node[second,label=left:$d_L$] at (0,0.7) {};
\node[first,label=above left:$a_L$] at (2.8,1.9) {};
\node[second,label=below left:$b_L$] at (2.8,0.7) {};
\node[second,label=above right:$y_L$] at (4.5,1.9) {};
\node[first,label=below right:$x_L$] at (4.5,0.7) {};
\node[second,label=above left:$y_R$] at (7.5,1.9) {};
\node[first,label=below left:$x_R$] at (7.5,0.7) {};
\node[first,label=above right:$a_R$] at (9.2,1.9) {};
\node[second,label=below right:$b_R$] at (9.2,0.7) {};
\node[first,label=right:$c_R$] at (12,1.9) {};
\node[second,label=right:$d_R$] at (12,0.7) {};
\end{tikzpicture}
\caption{A schematic representation of the graph $X(S)$.}
\label{fig:X}
\end{figure}

If $S$ is nonempty, define $X(S)$ by adding the four edges $a_Ly_L$, $b_Lx_L$, $a_Ry_R$ and $b_Rx_R$ to the disjoint union of $D_L$, $S$ and $D_R$, as shown in Figure~\ref{fig:X}. If $S$ is empty, define $X(S)=X(\varnothing)$ by adding the edges $a_Lb_R$ and $b_La_R$ to the disjoint union of $D_L$ and $D_R$; this graph has $40$ vertices. Hence, $X(S)$ is a bipartite graph with partite sets $\mathcal A_X=\mathcal A_L\cup U\cup\mathcal A_R$ and $\mathcal B_X=\mathcal B_L\cup W\cup\mathcal B_R$ such that $|\mathcal A_X|=|\mathcal B_X|$. We call $c_L,d_L,c_R$ and $d_R$ the four terminals of $X(S)$. Note that $c_L,c_R \in \mathcal A_X$ and $d_L,d_R \in \mathcal B_X$.

For every graph $S$ satisfying the properties mentioned above, the following properties hold for the terminals of $X(S)$.

\begin{lemma}\label{lem:X}
The following statements hold for the graph $X(S)$:
\begin{enumerate}
\item[(i)] The only two distinct terminals that occur as endpoints of a Hamiltonian path are $\{c_L,d_R\}$ and $\{c_R,d_L\}$.
\item[(ii)] There do not exist two vertex-disjoint paths that together contain all vertices in $X(S)$ and have $c_L,d_L,c_R$ and $d_R$ as their four endpoints.
\item[(iii)] For every $z\in\mathcal A_X$, the graph $X(S)-z$ has a Hamiltonian $d_Ld_R$-path, and for every $z\in\mathcal B_X$, the graph $X(S)-z$ has a Hamiltonian $c_Lc_R$-path.
\end{enumerate}
\end{lemma}

\begin{proof}
We first prove (i). Consider a Hamiltonian $c_Lb_L$-path in $D_L$ and a Hamiltonian $a_Rd_R$-path in $D_R$ (which exist by Lemma~\ref{lem:D}(i)). If $S=\varnothing$, one obtains a Hamiltonian $c_Ld_R$-path by adding the edge $b_La_R$. Otherwise, one obtains a Hamiltonian $c_Ld_R$-path by adding the edges $b_Lx_L$, $y_Ra_R$ and a Hamiltonian $x_Ly_R$-path in $S$. Analogously, one can obtain a Hamiltonian $c_Rd_L$-path. Since $X(S)$ is a bipartite graph for which $|\mathcal A_X|=|\mathcal B_X|$, the endpoints of a Hamiltonian path must belong to different partite sets and therefore neither $\{c_L,c_R\}$ nor $\{d_L,d_R\}$ can be its endpoint set.

To prove (i), we still need to consider the endpoint sets $\{c_L,d_L\}$ and $\{c_R,d_R\}$. The proof of (ii) uses the same argument as the exclusion of these two endpoint sets, so we will treat them together.

Let $Q$ be a spanning subgraph of $X(S)$ such that $Q$ consists of either a Hamiltonian  $c_Ld_L$-path or of two vertex-disjoint paths with $c_L,d_L,c_R$ and $d_R$ as their four endpoints. Both $c_L$ and $d_L$ have degree one in $Q$, while all other vertices of $D_L$ have degree two in $Q$. Therefore, the sum of the degrees in $Q$ over all vertices of $D_L$ is even. Therefore, $Q$ must contain an even number of edges with one endpoint in $V(D_L)$ and another endpoint in $V(X(S)) \setminus V(D_L)$. The only two possible such edges are incident with $a_L$ and $b_L$.

If $Q$ contains neither of these two edges, the graph $Q[V(D_L)]$ would be a Hamiltonian $c_Ld_L$-path in $D_L$, which is impossible by Lemma~\ref{lem:D}(i). So $Q$ contains both of these two edges. This implies in turn that $Q[V(D_L)]$ consists of two vertex-disjoint paths that together contain all vertices of $V(D_L)$ and have $a_L,b_L,c_L$ and $d_L$ as their four endpoints. Lemma~\ref{lem:D}(ii) then implies that one of these is an $a_Lb_L$-path and the other one is a $c_Ld_L$-path. The $c_Ld_L$-path must then be a connected component of $Q$, since no vertices in this path are incident with an edge that leaves $D_L$. So $Q$ cannot be a Hamiltonian $c_Ld_L$-path in $X(S)$. Analogously, one can apply the same argument to $D_R$, which shows that a Hamiltonian $c_Rd_R$-path in $X(S)$ cannot exist and this completes the proof of (i).

If $Q$ consists of two vertex-disjoint paths with $c_L,d_L,c_R$ and $d_R$ as their four endpoints, our argument shows that one component is a $c_Ld_L$-path contained in $D_L$ and another component is a $c_Rd_R$-path contained in $D_R$. Neither of these two components contains $a_L$ or $b_L$, but then $Q$ does not contain all vertices of $X(S)$ and we are done with (ii).

Finally, we prove (iii). We first show that for every $z\in\mathcal A_X$, the graph $X(S)-z$ has a Hamiltonian $d_Ld_R$-path. We first treat the case where $z\in\mathcal A_L$. There is a Hamiltonian $d_Lb_L$-path in $D_L-z$ by Lemma~\ref{lem:D}(iii) and a Hamiltonian $a_Rd_R$-path in $D_R$ by Lemma~\ref{lem:D}(i). If $S=\varnothing$, adding the edge $b_La_R$ yields the required Hamiltonian $d_Ld_R$-path in $X(S)-z$. If $S \neq \varnothing$, adding the edges $b_Lx_L$, $y_Ra_R$ and a Hamiltonian $x_Ly_R$-path in $S$ yields the required Hamiltonian $d_Ld_R$-path in $X(S)-z$. With a very similar argument, we can deal with the case where $z\in\mathcal A_R$. We now treat the case where $z\in U$. Then we have $S \neq \varnothing$. The graph $S-z$ has a Hamiltonian $y_Ly_R$-path. We can add the edges $a_Ly_L$ and $y_Ra_R$ and a Hamiltonian $d_La_L$-path in $D_L$ and a Hamiltonian $a_Rd_R$-path in $D_R$ to obtain a Hamiltonian $d_Ld_R$-path in $X(S)-z$.

The second statement of (iii) follows by a completely analogous argument.
\end{proof}

\subsection{The Petersen graph construction}\label{subsec:Petersen}

We now combine five of the preceding graphs to obtain bipartite hypotraceable graphs. Let $P$ be the Petersen graph, using the labels as shown in Figure~\ref{fig:Petersen}. Consider the perfect matching $M=\{pp',qq',rr',ss',tt'\}$. The graph $P-M$ is the disjoint union of two cycles of length five. We orient their edges as $p\to q\to r\to s\to t\to p$ and $p'\to s'\to q'\to t'\to r'\to p'$. Let $\vec C$ denote the set of these ten directed edges.

\begin{center}
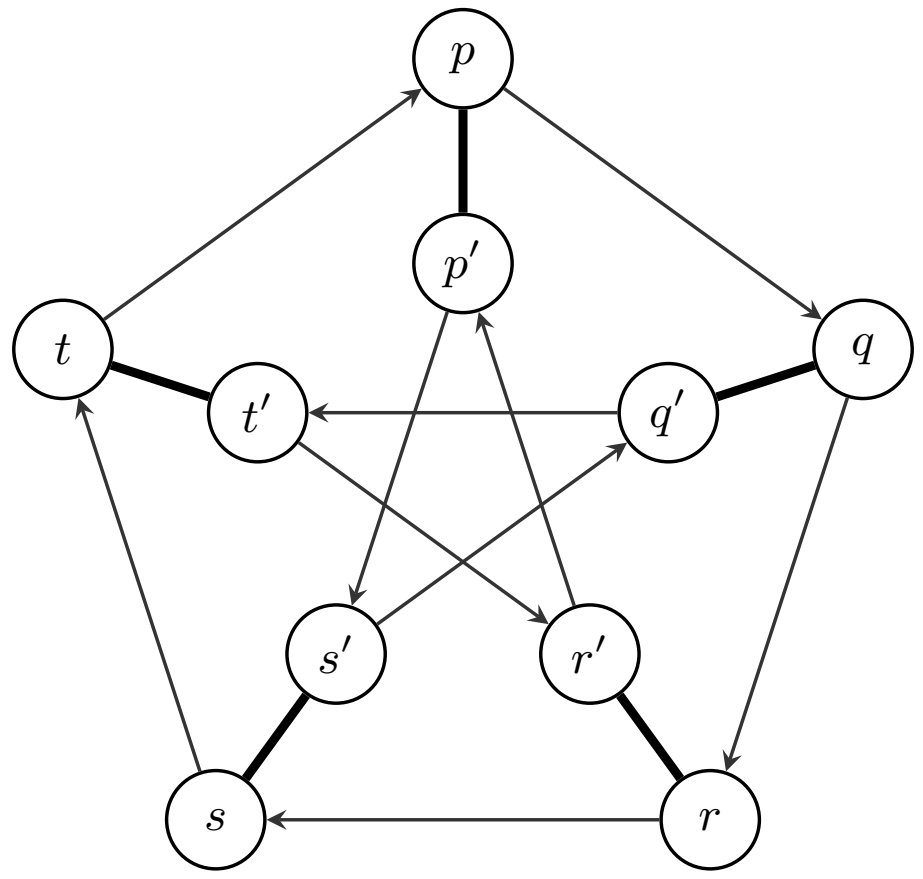

\begin{minipage}{\linewidth}
\centering
\begin{tikzpicture}[
 scale=1.05,>=stealth,
 vertex/.style={circle,draw=black,line width=0.65pt,fill=white,minimum size=6.5mm,inner sep=0pt,font=\small},
 cycle edge/.style={draw=black!80,line width=0.65pt,->},
 matching edge/.style={draw=black,line width=1.7pt}
]
\foreach \i/\name in {0/p,1/q,2/r,3/s,4/t}{
 \pgfmathsetmacro{\angle}{90-72*\i}
 \node[vertex] (outer\i) at (\angle:2.65cm) {$\name$};
 \node[vertex] (inner\i) at (\angle:1.36cm) {$\name'$};
}
\foreach \i in {0,...,4}{
 \pgfmathtruncatemacro{\j}{mod(\i+1,5)}
 \pgfmathtruncatemacro{\k}{mod(\i+3,5)}
 \draw[cycle edge] (outer\i)--(outer\j);
 \draw[cycle edge] (inner\i)--(inner\k);
 \draw[matching edge] (outer\i)--(inner\i);
}
\end{tikzpicture}
\captionof{figure}{The Petersen graph $P$ with the perfect matching $M$ in bold.}
\label{fig:Petersen}
\end{minipage}
\end{center}

The following observation can readily be seen.

\begin{obs}\label{obs:Petersen-matching}
The subgraph induced by the six endpoints of any three distinct edges of $M$ contains a cycle of length five.
\end{obs}

For each edge $e=uv \in M$, consider a graph $X(S_e)$ as in Lemma~\ref{lem:X} (for some valid choice of a graph $S_e$). The two copies of $D$ in the graph $X(S_e)$ are denoted by $D_u$ and $D_v$ and we refer to their four terminals as, respectively, $a_u,b_u,c_u,d_u$ and $a_v,b_v,c_v,d_v$. Therefore, the terminals of $X(S_e)$ are $c_u,d_u,c_v,d_v$. Note that the five graphs $X(S_e)$ do not need to be isomorphic.

Finally, let $G$ be the graph obtained by taking the disjoint union of these five graphs $X(S_e)$ and adding the edge $c_ud_v$ for each directed edge $u\to v$ in $\vec C$. These ten particular edges are called \textit{attachment edges}. For each vertex $u \in V(P)$, the outgoing edge at $u$ corresponds to an attachment edge incident with $c_u$, whereas the incoming edge at $u$ corresponds to an attachment edge incident with $d_u$. We are now ready to prove the main result of this paper.

\begin{theorem}\label{thm:main}
Every graph $G$ constructed in this way is bipartite and hypotraceable.
\end{theorem}

\begin{proof}
In each graph $X(S_e)$, the terminals $c_u,c_v$ belong to one partite set and the terminals $d_u,d_v$ belong to the other partite set. Therefore, $G$ is bipartite, since every attachment edge is between two vertices belonging to different partite sets.

Suppose for the sake of obtaining a contradiction that $G$ has a Hamiltonian path $Q$. Let $M_0$ consist of all edges $e\in M$ for which $X(S_e)$ contains neither endpoint of $Q$. We have $|M_0|\geq 3$. Let $P_0$ be the subgraph of $P$ that is induced by all endpoints of edges in $M_0$.

Consider an edge $e=uv \in M_0$. Every vertex of $X(S_e)$ must have degree two in $Q$. Therefore, an even number of attachment edges of $Q$ leaves $X(S_e)$. This number must be strictly larger than zero, because $Q$ is connected. There are only four attachment edges incident with $X(S_e)$. The graph $Q$ cannot contain all four attachment edges, because then the graph $Q[V(X(S_e))]$ would consist of two vertex-disjoint paths that together contain all vertices in $X(S_e)$ and have the four terminals as its endpoints. This is impossible by Lemma~\ref{lem:X}(ii). Hence, $Q$  contains exactly two of these attachment edges. Therefore, $Q[V(X(S_e))]$ is a Hamiltonian path in $X(S_e)$. The endpoints of this Hamiltonian path are either $c_u,d_v$ or $d_u,c_v$ by Lemma~\ref{lem:X}(i).

Therefore, for each vertex $w$ of $P_0$, exactly one of the incoming and outgoing edges of $\vec C$ corresponds to an attachment edge of $Q$. We give $w$ colour~$1$ if it is the outgoing edge, and colour~$0$ otherwise. By the preceding paragraph, the endpoints of every edge of $M_0$ have a different colour.  If $u\to v\in\vec C$ is such that both $u \in P_0$ and $v \in P_0$, then $u$ has colour~$1$ if and only if $c_ud_v\in E(Q)$, which in turn holds if and only if $v$ has colour~$0$. The endpoints of every edge of $P_0$ have a different colour, since this edge belongs either to $M_0$ or to $P-M$. Therefore, $P_0$ is bipartite. But since $|M_0| \geq 3$, $P_0$ contains a cycle of length five due to Observation~\ref{obs:Petersen-matching}. We reach a contradiction and therefore $G$ is not traceable.

Finally, we show that $G-z$ is traceable for every $z$ in $V(G)$. Since the five matching edges play symmetric roles in
the argument, we may relabel the vertices of the Petersen graph, together with their associated building blocks,
so that $z\in V(X(S_{pp'}))$. Consider the following two Hamiltonian paths of $P$: $ Q_{\mathrm{in}}=(s,s',q',q,r,r',p',p,t,t')$ and $ Q_{\mathrm{out}}=(r,r',t',t,s,s',p',p,q,q')$. Both of these Hamiltonian paths alternate between using one of the five edges of $M$ and an edge of $P-M$. The idea is to construct a Hamiltonian path in $G-z$ by replacing every edge of $Q_{\mathrm{in}}$ or $Q_{\mathrm{out}}$ by an appropriate path in $G$.

If $z$ is in the same partite set as $c_p$ and $c_{p'}$, we use $Q_{\mathrm{in}}$. The two edges of $P-M$ incident with the edge $pp'$ correspond to attachment edges incident with $d_p$ and $d_{p'}$. Therefore, Lemma~\ref{lem:X}(iii) yields the existence of a Hamiltonian $d_pd_{p'}$-path in $X(S_{pp'})-z$. If $z$ is in the other partite set, we use $Q_{\mathrm{out}}$. Now the attachment edges are incident with $c_p$ and $c_{p'}$ and Lemma~\ref{lem:X}(iii) yields the existence of a Hamiltonian $c_pc_{p'}$-path in $X(S_{pp'})-z$.

Now we replace every edge $e=uv \in M \setminus \{pp'\}$ of the chosen Hamiltonian path in $P$ by a Hamiltonian path in $X(S_e)$. If $e$ is not the first and not the last edge on the chosen Hamiltonian path in $P$, the required endpoints are one $c$-terminal and one $d$-terminal, belonging to different copies $D_u,D_v$. Therefore, Lemma~\ref{lem:X}(i) yields the necessary path. The first and the last edge on the chosen Hamiltonian path in $P$ both belong to $M$. 
For each of these two edges, let $u$ be its endpoint incident with the adjacent edge of $P-M$, and let $v$ be its other endpoint. If the corresponding attachment edge is incident with $c_u$, choose a Hamiltonian $c_ud_v$-path in $X(S_e)$. If it is incident with $d_u$, choose a Hamiltonian $d_uc_v$-path instead. Both choices are possible by Lemma~\ref{lem:X}(i). The endpoint in $D_v$ becomes an endpoint of the path we are constructing in $G-z$.

We now have one path covering each of the five graphs $X(S_e)$, with $z$ omitted from $X(S_{pp'})$. Join these paths using the four attachment edges corresponding to the edges of $P-M$ on the chosen Hamiltonian path of $P$. By construction, each attachment edge joins the appropriate endpoints of two consecutive paths. The resulting path contains every vertex of $G-z$ exactly once and is therefore a Hamiltonian path in $G-z$.
\end{proof}

\subsection{Consequences}\label{subsec:consequences}

Finally, we draw the reader's attention to several consequences as a result of the considerable freedom in our construction due to the choice of $X(S)$. If all five graphs in our construction are chosen as $X(\varnothing)$, we obtain a bipartite hypotraceable graph $G_0$ with $200$ vertices and $310$ edges. This is the smallest graph in our infinite family. If one of these five graphs is replaced by $X(S)$ (with $S \neq \varnothing$), we obtain a bipartite hypotraceable graph $G_S$ with $200+|V(S)|$ vertices and $312+|E(S)|$ edges.

For $r \geq 1$, let $L_r$ be a bipartite graph with partite sets $\{x_1,\ldots,x_r\}$ and $\{y_1,\ldots,y_r\}$ such that $x_iy_j\in E(L_r)$ if and only if $|i-j| \leq 1$. Choose $x_L=x_1$, $x_R=x_r$, $y_L=y_1$ and $y_R=y_r$. The graph $L_r$ satisfies conditions~(\ref{cond:S-paths}) and~(\ref{cond:S-deletion}): the paths $(x_1,y_1,\ldots,x_r,y_r)$ and $(y_1,x_1,\ldots,y_r,x_r)$ establish the first condition, while $(y_1,x_1,\ldots,y_{j-1},x_{j-1},y_j,x_{j+1},y_{j+1},\ldots,x_r,y_r)$ is a Hamiltonian $y_1y_r$-path in $L_r-x_j$ (when $j=1$ or $j=r$, the initial or final segment is omitted, respectively). The other paths for the second condition are obtained by exchanging $x_i$ and $y_i$. All these paths of course still exist whenever an arbitrary set of edges is added and therefore every bipartite graph $S$ for which $L_r$ is a subgraph of $S$ and $S$ is a subgraph of $K_{r,r}$ yields a valid choice for the graph $S$ required in Lemma~\ref{lem:X}. This has the following interesting consequences.

\begin{theorem}\label{thm:edges}
The following statements hold:
\begin{enumerate}
\item[(i)] For every even $n \geq 200$, there is a bipartite hypotraceable graph of order $n$ with exactly $20$ vertices of degree four and all remaining vertices of degree three. For even $n \geq 202$, every size between $3n/2+10$ and $(n-200)^2/4+312$ is attained.
\item[(ii)]\label{thm:asymptotics} As $n$ tends to infinity through even integers, the maximum size of a bipartite hypotraceable graph of order $n$ is $\left(1/4+o(1)\right)n^2$, and the number $h(n)$ of nonisomorphic such graphs satisfies
\[
 h(n)\geq2^{(1/4+o(1))n^2}.
\]
\item[(iii)]\label{thm:induced} Every bipartite graph on $s \geq 1$ vertices is an induced subgraph of a bipartite hypotraceable graph on $2s+202$ vertices.
\end{enumerate}
\end{theorem}

\begin{proof}
For (i), use $G_0$ when $n=200$, and otherwise put $r=(n-200)/2$ and use $S=L_r$ for the graph $G_S$. The $4 \cdot 5=20$ terminals have degree four, while all other vertices have degree three. Since $|E(L_r)|=3r-2$, adding any number of its missing edges between vertices in distinct partite sets gives every size from $312+3r-2$ to $312+r^2$.

For (ii), taking $S=K_{r,r}$, where $r=(n-200)/2$, gives graphs with $312+r^2=\left(1/4+o(1)\right)n^2$ edges. The density assertion follows from the upper bound $n^2/4$ for the size of any bipartite graph of order $n$.

For the counting assertion, independently including or omitting each of the $r^2-3r+2$ edges outside $L_r$ gives $2^{r^2-3r+2}$ distinct labelled graphs $G_S$. Each isomorphism class contains at most $n!$ labelled graphs, so
\[
 h(n)\geq\frac{2^{r^2-3r+2}}{n!}
       =2^{(1/4+o(1))n^2},
\]
since $r=n/2-100$ and $\log_2(n!)=o(n^2)$.

For (iii), let the prescribed graph have partite sets of sizes $\alpha,\beta$, with $\alpha+\beta=s$, and take $r=s+1$. The sets $A=\{x_1,\ldots,x_\alpha\}$ and $B=\{y_{\alpha+2},\ldots,y_{s+1}\}$ (with either set empty when its size is zero) together form an independent set in $L_r$, since $j-i\geq2$ for every $x_i\in A$ and $y_j\in B$. Add the prescribed edges between $A$ and $B$ to obtain a graph $S$ (such that $L_r$ is a subgraph of $S$ and $S$ is a subgraph of $K_{r,r}$). Thus $S$ yields a valid choice for the graph $S$ required in Lemma~\ref{lem:X} and $G_S$ contains the prescribed graph as an induced subgraph, where $|V(G_S)|=200+2r=2s+202$.
\end{proof}

We note that Chv\'atal~\cite{Chvatal73} asked in 1973 whether every graph occurs as an induced subgraph of some hypohamiltonian graph. This question was answered affirmatively by Zamfirescu and Zamfirescu~\cite{ZZ18}. Hence, Theorem~\ref{thm:edges}(iii) yields
an analogous result for bipartite hypotraceable graphs.

\section{Open questions}

The search for smaller counterexamples to a conjecture can lead to further structural insights. Tutte~\cite{Tutte71} famously conjectured that every $3$-connected cubic bipartite graph is Hamiltonian. Horton's $96$-vertex counterexample, recorded by Bondy and Murty~\cite[p.~240]{BondyMurty76}, was followed by smaller examples, including the $54$-vertex graph of Ellingham and Horton~\cite{EH83}. Georges~\cite{Georges89} and Kelmans~\cite{Kelmans88} independently constructed a counterexample on $50$ vertices. Brinkmann, Goedgebeur and McKay~\cite{BGM22} eventually proved that this Georges--Kelmans graph is a smallest counterexample; they also constructed an infinite family of counterexamples with cyclic edge-connectivity at least five. This progression motivates the corresponding search for smaller bipartite hypotraceable graphs.

\begin{question}
    What is the smallest order of a bipartite hypotraceable graph?
\end{question}

The construction above yields a bipartite hypotraceable graph on $200$ vertices. A small modification yields a bipartite hypotraceable graph on $150$ vertices, which we describe here without proof. The verification of this smaller example is a bit more cumbersome and we have chosen to instead develop the preceding construction for the sake of clarity of its proof rather than to minimise the order.

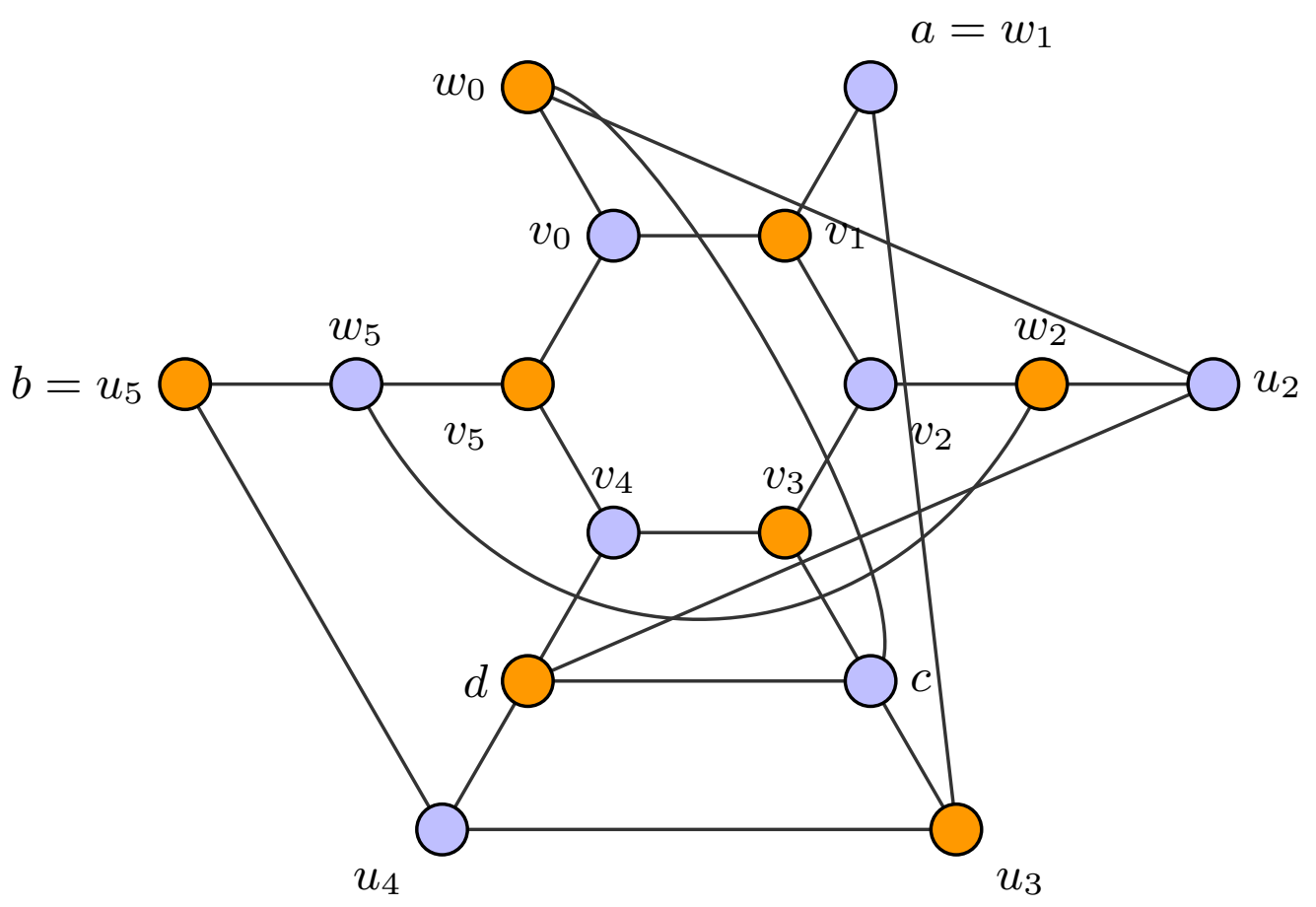
\begin{figure}[htbp]
\centering
\begin{tikzpicture}[
  x=1.15cm,y=1.15cm,
  edge/.style={draw=black!80,line width=0.65pt},
  vertex/.style={circle,draw=black,line width=0.65pt,inner sep=0pt,minimum size=3.4mm},
  first/.style={vertex,fill=blue!25},
  second/.style={vertex,fill=orange!80!yellow}
]
\coordinate (u2) at (3,0);
\coordinate (u3) at (1.5,-2.598076);
\coordinate (u4) at (-1.5,-2.598076);
\coordinate (u5) at (-3,0);
\coordinate (w0) at (-1,1.732051);
\coordinate (w1) at (1,1.732051);
\coordinate (w2) at (2,0);
\coordinate (c) at (1,-1.732051);
\coordinate (d) at (-1,-1.732051);
\coordinate (w5) at (-2,0);
\coordinate (v0) at (-0.5,0.866025);
\coordinate (v1) at (0.5,0.866025);
\coordinate (v2) at (1,0);
\coordinate (v3) at (0.5,-0.866025);
\coordinate (v4) at (-0.5,-0.866025);
\coordinate (v5) at (-1,0);

\foreach \u/\v in {
  w1/v1, v2/v1, v2/v3, v2/w2, u2/w2,
  v0/v5, v0/w0, c/v3, c/u3,
  v4/v3, v4/d, v4/v5,
  u4/u3, u4/d, u4/u5,
  w5/v5, w5/u5,
  w0/u2, w1/u3, d/u2, v0/v1, c/d} {
  \draw[edge] (\u) -- (\v);
}
\draw[edge] (w5) .. controls (-1.1,-1.83) and (1.1,-1.83) .. (w2);
\draw[edge] (w0) .. controls (-0.45,2.01) and (1.51,-1.40) .. (c);

\foreach \x in {u2,u4,w1,w5,v0,v2,v4,c}
  \node[first] at (\x) {};
\foreach \x in {u3,u5,w0,w2,d,v1,v3,v5}
  \node[second] at (\x) {};

\node[right=4pt,font=\small] at (u2) {$u_2$};
\node[below right=4pt,font=\small] at (u3) {$u_3$};
\node[below left=4pt,font=\small] at (u4) {$u_4$};
\node[left=4pt,font=\small] at (u5) {$b=u_5$};
\node[left=4pt,font=\small] at (w0) {$w_0$};
\node[above right=4pt,font=\small] at (w1) {$a=w_1$};
\node[above=4pt,font=\small] at (w2) {$w_2$};
\node[right=4pt,font=\small] at (c) {$c$};
\node[left=4pt,font=\small] at (d) {$d$};
\node[above=4pt,font=\small] at (w5) {$w_5$};
\node[left=4pt,font=\small] at (v0) {$v_0$};
\node[right=4pt,font=\small] at (v1) {$v_1$};
\node[below right=4pt,font=\small] at (v2) {$v_2$};
\node[above=4pt,font=\small] at (v3) {$v_3$};
\node[above=4pt,font=\small] at (v4) {$v_4$};
\node[below left=4pt,font=\small] at (v5) {$v_5$};
\end{tikzpicture}
\caption{The $16$-vertex graph $K$, with terminals $a=w_1$, $b=u_5$, $c$ and $d$.}
\label{fig:K}
\end{figure}

Starting from the graph $B$ in Figure~\ref{fig:B}, delete the vertices $u_0$ and $u_1$, replace the edges $w_0w_1$ and $u_2u_3$ by $w_0u_2$ and $w_1u_3$, and add the edges $cw_0$, $du_2$ and $v_0v_1$. Denote the resulting $16$-vertex graph by $K$, with terminals $a=w_1$, $b=u_5$, $c$ and $d$ (see Figure~\ref{fig:K}).

Take two disjoint copies $K_L$ and $K_R$, and identify $a_L$ with $a_R$ and $b_L$ with $b_R$. The resulting $30$-vertex graph $F$, with terminals $c_L,d_L,c_R,d_R$ also satisfies the three conclusions of Lemma~\ref{lem:X}. Therefore, we can replace each of the five graphs $X(S_e)$ in the Petersen construction by a copy of $F$, yielding a bipartite hypotraceable graph on $150$ vertices.

Our second question is motivated by a paper of Thomassen~\cite{Thomassen76}. He constructed infinite families of planar hypotraceable graphs in $1976$, while the present paper establishes the existence of bipartite hypotraceable graphs. It is natural to ask whether these two properties can be combined.

\begin{question}
    Does there exist a planar bipartite hypotraceable graph?
\end{question}

\section*{Use of AI}

The starting point for this work was a bipartite hypotraceable graph obtained over several days of interaction with the AI model Fable 5.2. However, this output provided neither a structural description of the graph nor a conceptual explanation of these properties. The initial proof that this graph was hypotraceable was entirely computer-assisted: a certificate of unsatisfiability for a SAT encoding established non-traceability, while a list of Hamiltonian paths established traceability for each vertex-deleted subgraph. The author subsequently identified the special role of the 18-vertex Ellingham–Horton graph and the Petersen graph, developed alternative constructions of bipartite hypotraceable graphs, and introduced the auxiliary bipartite graph $S$ in the graph $X(S)$ that led to the consequences mentioned in this paper. The author also developed, rewrote and simplified the proofs and exposition. AI tools assisted at later stages with proofreading and drafting figures.

\end{document}